\documentclass{amsart}
\usepackage{amsmath}
\usepackage{amsthm}
\usepackage{amssymb}

\theoremstyle{plain}
\newtheorem{thm}{Theorem}[section]
\newtheorem{prop}[thm]{Proposition}

\newtheorem{lem}[thm]{Lemma}
\theoremstyle{remark}

\theoremstyle{definition}
\newtheorem{exmp}[thm]{Example}
\newtheorem{defn}[thm]{Definition}
\numberwithin{equation}{section}
\DeclareMathOperator*{\slim}{s-lim}

\begin{document}
\title{Prediction of fractional processes with long-range dependence}
\author{Akihiko Inoue}
\author{Vo V. Anh}

\address{Department of Mathematics, Hiroshima University, 
Higashi-Hiroshima 739-8526, Japan}
\email{inoue100@hiroshima-u.ac.jp}
\address{School of Mathematical Sciences, 
Queensland University of Technology, 
GPO Box 2434, Brisbane, Queensland 4001, 
Australia}
\email{v.anh@qut.edu.au}
\date{9 November 2011}

\begin{abstract}
We introduce a class of Gaussian processes with stationary increments which
exhibit long-range dependence. The class includes fractional Brownian motion
with Hurst parameter $H>1/2$ as a typical example. We establish infinite and
finite past prediction formulas for the processes in which the predictor
coefficients are given explicitly in terms of the MA$(\infty )$ and AR$%
(\infty )$ coefficients.
\end{abstract}

\maketitle

\renewcommand{\thefootnote}{\fnsymbol{footnote}} \footnote[0]{%
2000\textit{\ Mathematics Subject Classification} Primary 60G25; Secondary
60G15.}


\section{Introduction}
\label{sec:1}

Let $(X(t):t\in \mathbf{R})$ be a centered Gaussian process with stationary
increments, defined on a probability space $(\Omega ,\mathcal{F},P)$, that
admits the \textit{moving-average} representation 
\begin{equation}
X(t)=\int_{-\infty }^{\infty }\left\{ g(t-s)-g(-s)\right\} dW(s),\qquad t\in 
\mathbf{R},  \label{eq:1.1}
\end{equation}%
where $(W(t):t\in \mathbf{R})$ is a Brownian motion, and $g(t)$ is a
function of the form 
\begin{align}
& g(t)=\int_{0}^{t}c(s)ds,\qquad t\in \mathbf{R},  \label{eq:1.2} \\
& c(t):=I_{(0,\infty )}(t)\int_{0}^{\infty }e^{-ts}\nu (ds),\qquad t\in 
\mathbf{R},  \label{eq:1.3}
\end{align}%
with some Borel measure $\nu $ on $(0,\infty )$ satisfying 
\begin{equation}
\int_{0}^{\infty }\frac{1}{1+s}\nu (ds)<\infty .  \label{eq:1.4}
\end{equation}%
We will also assume some extra conditions such as 
\begin{align}
& \lim_{t\rightarrow 0+}c(t)=\infty ,  \label{eq:1.5} \\
& g(t)\sim t^{H-(1/2)}\ell (t)\cdot \frac{1}{\Gamma (\frac{1}{2}+H)},\qquad
t\rightarrow \infty ,  \label{eq:1.6}
\end{align}%
where $\ell (t)$ is a slowly varying function at infinity and $H$ is a
constant such that 
\begin{equation}
1/2<H<1.  \label{eq:1.7}
\end{equation}%
In (\ref{eq:1.6}), and throughout the paper, $a(t)\sim b(t)$ as $%
t\rightarrow \infty $ means $\lim_{t\rightarrow \infty }a(t)/b(t)=1$. We
call $c(t)$ (as well as $g(t)$) the \textit{MA$(\infty)$ coefficient} of 
$(X(t))$. 
We remark that, in the prediction formulas for $(X(t))$ which we consider in this paper, 
$c(t)$ becomes more relevant than $g(t)$.

A typical example of $\nu $ is 
\begin{equation}
\nu (ds)=\frac{\sin \{\pi (H-\frac{1}{2})\}}{\pi }s^{(1/2)-H}ds\quad 
\mbox{on
$(0,\infty)$}  \label{eq:1.8}
\end{equation}%
with (\ref{eq:1.7}). For this $\nu $, $g(t)$ becomes 
\begin{equation}
g(t)=I_{(0,\infty )}(t)t^{H-(1/2)}\frac{1}{\Gamma (\frac{1}{2}+H)},\qquad t\in 
\mathbf{R},  \label{eq:1.9}
\end{equation}%
and $(X(t))$ reduces to \textit{fractional Brownian motion\/} $(B_{H}(t))$
with \textit{Hurst parameter} $H$ (see Example \ref{exmp:2.3} below).
Fractional Brownian motion, abbreviated fBm, was introduced by Kolmogorov 
\cite{K}. For $1/2<H<1,$ fBm has both \textit{self-similarity} and \textit{%
long-range dependence }(Samorodnitsky and Taqqu \cite{ST}\textit{\/}), and
plays an important role in various fields such as network traffic (see,
e.g., Mikosch et al.\/ \cite{MRRS}) and finance (see, e.g., Hu et al.\ \cite%
{HOS}); see also Taqqu \cite{T} and other papers in the same volume. Because
of its importance, stochastic calculus for fBm has been developed by many
authors; see, e.g., Decreusefond and \"{U}st\"{u}nel \cite{DU}, and Nualart 
\cite{N}. Grecksch and Anh \cite{GA99} introduced Hilbert space-valued fBm
and the corresponding stochastic calculus. Duncan et al.\ \cite{DMP02}
and Tindel et al.\ \cite{TTV03} studied stochastic evolution equations
with fBm in Hilbert spaces. Other important examples of $(X(t))$ are the
processes with long-range dependence which, unlike fBm, have two different
indices $H_{0}$ and $H$ describing the local properties (path properties)
and long-time behavior of $(X(t))$, respectively (see Example \ref{exmp:2.4}
below).

Let $t_0$, $t_1$ and $T$ be real constants such that 
\begin{equation}
-\infty<-t_0\le 0\le t_1< T<\infty,\qquad -t_0<t_1.  \label{eq:1.10}
\end{equation}
For $I=(-\infty,t_1]$ or $[-t_0,t_1]$, we write $P_IX(T)$ for the predictor
of the future value $X(T)$ based on the observable $(X(s): s\in I)$ (see
\S 3 below). One of the fundamental prediction problems for $(X(t))$ is
to express $P_IX(T)$ using the segment $(X(s):s\in I)$ and some
deterministic quantities. Another is to express the variance of the
prediction error $P_I^{\bot}X(T):=X(T)-P_IX(T)$. Results of this type become
important tools in the analysis of non-Markovian processes and systems
modulated by them (see, e.g., Norros et al.\ \cite{NVV}, Anh et al.\ \cite%
{AIK}, Inoue et al.\ \cite{INA} and Inoue and Nakano \cite{IN}). One of our
main purposes here is to derive such results for $(X(t))$.

We establish the following infinite and finite past prediction formulas for $%
(X(t))$ (see Theorems \ref{thm:3.8} and \ref{thm:4.12} below): 
\begin{align}
& P_{(-\infty ,t_{1}]}X(T)=X(t_{1})+\int_{-\infty }^{t_{1}}\left\{
\int_{0}^{T-t_{1}}b(t_{1}-s,\tau )d\tau \right\} dX(s),  \label{eq:1.11} \\
& P_{[-t_{0},t_{1}]}X(T)=X(t_{1})+\int_{-t_{0}}^{t_{1}}\left\{
\int_{0}^{T-t_{1}}h(s+t_{0},u)du\right\} dX(s).  \label{eq:1.12}
\end{align}%
The significance of (\ref{eq:1.11}) and (\ref{eq:1.12}) is that the
predictor coefficients $b(t,s)$ and $h(t,s)$ are given explicitly in terms
of the MA$(\infty)$ coefficient $c(t)$ and AR$(\infty )$ coefficient $a(t)$, 
to be defined in \S \ref{subsec:3.1}, of $(X(t))$. 
The integral of $a(t)$ is in fact the coefficient of an AR$(\infty)$-type equation 
describing $(X(t))$ (see \S 5). 
We will find that $a(t)$ has a nice integral
representation similar to (\ref{eq:1.3}) (see (\ref{eq:3.3}) below). It
turns out that the existence of such a nice AR$(\infty )$ coefficient, in
addition to the nice MA$(\infty )$ coefficient, is a key to the solution to
the prediction problems above.

For fBm with $1/2<H<1$, the predictor coefficients $b(t,s)$ and $h(t,s)$ are
given in Gripenberg and Norros \cite{GN}. See \cite{NVV} and \cite{NP} for
different proofs. Fractional Brownian motion has a variety of nice
properties, and the methods of proof of \cite{GN,NVV,NP} naturally rely on
such special properties of fBm, hence are not applicable to $(X(t))$. The
method of this paper is based on the \textit{alternating projections to the
past and future\/} (see \S 4.1 below).
As for fBm with $0<H<1/2$, its
infinite and finite past prediction formulas also exist, and are due to
Yaglom \cite{Y} and Nuzman and Poor \cite{NP}, respectively (see also Anh
and Inoue \cite{AI1}).

In Inoue and Anh \cite{IA}, a class of processes $(\tilde{X}(t))$ 
of the same form
\begin{equation}
\tilde{X}(t)=\int_{-\infty }^{\infty }\left\{ \tilde{c}(t-s)-\tilde{c}(-s)\right\} dW(s),\qquad t\in 
\mathbf{R},  \label{eq:1.13}
\end{equation}
as (\ref{eq:1.1}) are introduced. 
Unlike $g(t)$ in (\ref{eq:1.1}), however, the kernel $\tilde{c}(t)$ itself is assumed to be of the form
\begin{equation}
\tilde{c}(t)=I_{(0,\infty )}(t)\int_{0}^{\infty }e^{-ts}\tilde{\nu} (ds),\qquad t\in 
\mathbf{R},  \label{eq:1.14}
\end{equation}
with a Borel measure $\tilde{\nu}$ on $(0,\infty)$ satisfying some suitable conditions. 
This class of $(\tilde{X}(t))$ includes fBm with $H\in (0, 1/2)$ as a typical example. 
Notice that $\tilde{c}(t)$ in (\ref{eq:1.14}) (resp., $g(t)$ in (\ref{eq:1.1})) is decreasing 
(resp., increasing) on $(0,\infty)$ as $t^{H-(1/2)}$ 
with $H\in (0,1/2)$ (resp., (1/2,1)) is. 
In \cite{IA}, prediction formulas for $(\tilde{X}(t))$ are proved, extending 
the results for fBm with $H\in (0,1/2)$ stated above. 
These prediction formulas for  $(\tilde{X}(t))$, including those for fBm with $H\in (0,1/2)$, 
have different forms from (\ref{eq:1.11}) and (\ref{eq:1.12}), 
in that no stochastic integrals appear there.

We provide the basic properties and examples of $(X(t))$ in \S 2. We
consider the infinite and finite past prediction problems for $(X(t))$ in
\S 3 and \S 4, respectively. 
Finally in \S 5, we remark on the AR$(\infty)$-type equations 
describing $(X(t))$ and $(\tilde{X}(t))$. 


\section{Basic properties and examples}
\label{sec:2}

In this section, we assume (\ref{eq:1.2})--(\ref{eq:1.4}) and 
\begin{equation}
\int_{1}^{\infty}c(t)^2dt<\infty.  \label{eq:2.1}
\end{equation}
Then, as in \cite[Lemma 2.1]{IA}, we have $\int_{-\infty}^{\infty}\vert
g(t-s)-g(-s)\vert^2ds<\infty$ for $t\in\mathbf{R}$. Therefore, for a
one-dimensional standard Brownian motion $(W(t):t\in\mathbf{R})$ with $%
W(0)=0 $, we may define the centered stationary-increment Gaussian process $%
(X(t):t\in\mathbf{R})$ by (\ref{eq:1.1}).

For $s>0$ and $t\in \mathbf{R}$, we put $\Delta _{s}X(t):=X(t+s)-X(t)$.
Then, by definition, $(\Delta _{s}X(t):t\in \mathbf{R})$ is a stationary
process.

\begin{lem}
\label{lem:2.1} Let $s\in (0,\infty)$. We assume\/ $(\ref{eq:1.6})$ and\/ $(%
\ref{eq:1.7})$. Then 
\begin{equation*}
E\left[\Delta_sX(t)\cdot\Delta_sX(0)\right] \sim t^{2H-2}\ell(t)^2\cdot 
\frac{s^2\Gamma(2-2H)\sin\{(H-\frac{1}{2})\pi\}}{\pi}, \qquad t\to\infty.
\end{equation*}
\end{lem}

Since $-1<2H-2<0$ in Lemma \ref{lem:2.1}, we see from this lemma that $%
(\Delta_sX(t))$, whence $(X(t))$, has long-range dependence.

We put $\sigma(t):=E[\vert X(t+s)-X(s)\vert^2]^{1/2}$ for $t\ge 0$ and $s\in%
\mathbf{R}$.

\begin{lem}
\label{lem:2.2} Let $H_{0}\in (1/2,1)$ and $\ell _{0}(\cdot )$ a slowly
varying function at infinity. We assume 
\begin{equation}
g(t)\sim t^{H_{0}-(1/2)}\ell_{0}(1/t)\cdot \frac{1}{\Gamma (\frac{1}{2}%
+H_{0})},\qquad t\rightarrow 0+.  \label{eq:2.2}
\end{equation}%
Then 
\begin{equation*}
\sigma (t)\sim t^{H_{0}}\ell_{0}(1/t)\sqrt{v(H_{0})},\qquad t\rightarrow 0+,
\end{equation*}%
where $v(H_{0}):=\Gamma (2-2H_{0})\cos (\pi H_{0})/\{\pi H_{0}(1-2H_{0})\}$.
In particular, we have 
\begin{equation*}
H_{0}=\sup \{\beta :\sigma (t)=o(t^{\beta }),\quad t\rightarrow 0+\}=\inf
\{\beta :t^{\beta }=o(\sigma (t)),\quad t\rightarrow 0+\}.
\end{equation*}
\end{lem}

From Lemma \ref{lem:2.2}, we see that the index $H_0$ describes the path
properties of $(X(t))$ (see Adler \cite[\S 8.4]{A}).

By the monotone density theorem (cf.\ Bingham et al.\ \cite[Theorem 1.7.5]%
{BGT}), (\ref{eq:1.6}) with (\ref{eq:1.7}) implies 
\begin{equation}
c(t)\sim t^{H-(3/2)}\ell(t)\cdot\frac{1}{\Gamma(H-\frac{1}{2})}, \qquad
t\to\infty.  \label{eq:2.3}
\end{equation}
Similarly, (\ref{eq:2.2}) implies 
\begin{equation}
c(t)\sim t^{H_0-(3/2)}\ell_0(1/t) \cdot\frac{1}{\Gamma(H_0-\frac{1}{2})}.
\qquad t\to 0+.  \label{eq:2.4}
\end{equation}
Lemmas \ref{lem:2.1} and \ref{lem:2.2} follow from (\ref{eq:2.3}) and (\ref%
{eq:2.4}), respectively, by standard arguments. However, since we do not use
these results, we omit the details.

\begin{exmp}
\label{exmp:2.3} For $H\in (1/2,1)$, let $\nu $ be as in (\ref{eq:1.8}).
Then we have (\ref{eq:1.9}); and so all the conditions above are satisfied.
The resulting process $(X(t))$ is fBm $(B_{H}(t))$: 
\begin{equation}
B_{H}(t)=\frac{1}{\Gamma (\frac{1}{2}+H)}\int_{-\infty }^{\infty }\left\{
((t-s)_{+})^{H-(1/2)}-((-s)_{+})^{H-(1/2)}\right\} dW(s),  \label{eq:2.5}
\end{equation}%
where $(x)_{+}:=\max (0,x)$ for $x\in \mathbf{R}$. The representation (\ref{eq:2.5})
of fBm is due to the pioneering work of Mandelbrot and Van Ness \cite{MV}.%
\textrm{\ }
\end{exmp}

\begin{exmp}
\label{exmp:2.4} Let $f(\cdot)$ be a nonnegative, locally integrable
function on $(0,\infty)$. For $H_0, H\in (1/2,1)$ and slowly varying
functions $\ell_0(\cdot)$ and $\ell(\cdot)$ at infinity, we assume 
\begin{align*}
&f(s)\sim \frac{\sin\{\pi (H_0-\frac{1}{2})\}}{\pi}s^{(1/2)-H}\ell(1/s),
\qquad s\to 0+, \\
&f(s)\sim \frac{\sin\{\pi (H_0-\frac{1}{2})\}}{\pi}s^{(1/2)-H_0}\ell_0(s),
\qquad s\to \infty.
\end{align*}
Let $\nu(ds)=f(s)ds$. Then, by Abelian theorems for Laplace transforms (cf.~%
\cite[\S 1.7]{BGT}), we have (\ref{eq:2.3}), whence (\ref{eq:1.6}).
Similarly, we have (\ref{eq:2.4}), whence (\ref{eq:2.2}). Thus all the
conditions above are satisfied. As we have seen above, the indices $H_0$ and 
$H$ describe the path properties and long-time behavior of $(X(t))$,
respectively.
\end{exmp}


\section{Infinite past prediction problems}
\label{sec:3}

In this section, we assume (\ref{eq:1.1})--(\ref{eq:1.5}), (\ref{eq:2.1})
and 
\begin{equation}
\lim_{t\to\infty}g(t)=\infty.  \label{eq:3.1}
\end{equation}
Notice that, for the processes $(X(t))$ in Examples \ref{exmp:2.3} and \ref%
{exmp:2.4}, all these conditions are satisfied. We also assume (\ref{eq:1.10}).

We write $M(X)$ for the real Hilbert space spanned by $\{X(t):t\in\mathbf{R}\}$
in $L^2(\Omega,\mathcal{F},P)$, and $\Vert \cdot\Vert$ for its norm. Let $I$
be a closed interval of $\mathbf{R}$ such as $[-t_0,t_1]$, $(-\infty,t_1]$,
and $[-t_0,\infty)$. Let $M_{I}(X)$ be the closed subspace of $M(X)$ spanned
by $\{X(t):t\in I\}$. We write $P_{I}$ for the orthogonal projection operator
from $M(X)$ to $M_{I}(X)$, and $P^{\bot}_{I}$ for its orthogonal complement: 
$P^{\bot}_{I}Z=Z-P_{I}Z$ for $Z\in M(X)$. Note that, since $(X(t))$ is a
Gaussian process, we have 
$P_{I}Z=E[Z \vert \sigma(X(s): s\in I)]$.

\subsection{MA and AR coefficients}\label{subsec:3.1}

The conditions (\ref{eq:1.5}) and (\ref{eq:3.1}) imply $\nu (0,\infty
)=\infty $ and $\int_{0}^{\infty }s^{-1}\nu (ds)=\infty $, respectively.
Therefore, by \cite[Theorem 3.2]{IA}, there exists a unique Borel measure $%
\mu $ on $(0,\infty )$ satisfying 
\begin{equation*}
\int_{0}^{\infty }\frac{1}{1+s}\mu (ds)<\infty ,\quad \mu (0,\infty )=\infty
,\quad \int_{0}^{\infty }\frac{1}{s}\mu (ds)=\infty
\end{equation*}%
and 
\begin{equation}
-iz\left\{ \int_{0}^{\infty }e^{izt}c(t)dt\right\} \left\{ \int_{0}^{\infty
}e^{izt}\alpha (t)dt\right\} =1,\qquad \Im z>0,  \label{eq:3.2}
\end{equation}%
with 
\begin{equation*}
\alpha (t):=\int_{0}^{\infty }e^{-st}\mu (ds),\qquad t>0.
\end{equation*}%
We define
\begin{equation}
a(t):=-\frac{d\alpha }{dt}(t)=\int_{0}^{\infty }e^{-st}s\mu (ds),\qquad t>0.
\label{eq:3.3}
\end{equation}
We call $a(t)$ (as well as $\alpha(t)$) the \textit{AR$(\infty)$ coefficient} of 
$(X(t))$ (see \S \ref{sec:5} for background). 
We define the positive kernel $b(t,s)$ by 
\begin{equation*}
b(t,s):=\int_{0}^{s}c(u)a(t+s-u)du,\qquad t,s>0.
\end{equation*}%
Then, by \cite[Lemma 3.4]{IA}, the following equalities hold: 
\begin{align}
& \int_{0}^{\infty }b(t,s)dt=1,\qquad s>0,  \label{eq:3.4} \\
& c(t+s)=\int_{0}^{t}c(t-u)b(u,s)du,\qquad t,s>0.  \label{eq:3.5}
\end{align}

\subsection{Stochastic integrals}\label{subsec:3.2}

Let $I$ be a closed interval of $\mathbf{R}$. We define 
\begin{equation*}
\mathcal{H}_{I}(X):=\left\{ f\ :\ 
\begin{aligned} &\mbox{$f$ is a real-valued measurable function on $I$ 
such}\\ &\mbox{that $\int_{-\infty}^{\infty}\left\{\int_{I}
\vert f(u)\vert c(u-s)du\right\}^2ds<\infty$.} \end{aligned}\right\} .
\end{equation*}%
This is the class of functions $f$ for which we can define the stochastic
integral $\int_{I}f(s)dX(s)$. 
We notice that, by Lemma \ref{lem:5.2} below, the function $c(t)$, 
whence $\mathcal{H}_{I}(X)$, is uniquely determined by $(X(t))$. 
We define a subclass $\mathcal{H}_{I}^{0}$ of $%
\mathcal{H}_{I}(X)$ by 
\begin{equation*}
\mathcal{H}_{I}^{0}:=\left\{
\sum\nolimits_{k=1}^{m}a_{k}I_{(t_{k-1},t_{k}]}(s):\begin{aligned}
&\mbox{$m\in \mathbf{N}$, $-\infty<t_{0}<t_{1}<\dots<t_{m}<\infty$} \\
&\mbox{with $(t_0, t_m]\subset I$, $a_{k}\in \mathbf{R}\ (k=1,\dots,m)$}
\end{aligned}\right\} .
\end{equation*}%
Each member of $f\in \mathcal{H}_{I}^{0}$ is a \textit{simple function\/} on $I$.

\begin{defn}
\label{def:3.1} For $f=\sum\nolimits_{k=1}^{m}a_{k}I_{(t_{k-1},t_{k}]}\in
\mathcal{H}_{I}^{0}$, we define 
\begin{equation*}
\int_{I}f(s)dX(s):=\sum_{k=1}^{m}a_{k}\left\{ X(t_{k})-X(t_{k-1})\right\}.
\end{equation*}
\end{defn}

We see that $\int_{I} f(s)dX(s)\in M_I(X)$ for $f\in \mathcal{H}_I^0$.

\begin{prop}
\label{prop:3.2} For $f\in\mathcal{H}_I^0$, we have 
\begin{equation}
\int_{I} f(s)dX(s)=\int_{-\infty}^{\infty} \left\{\int_{I}
f(u)c(u-s)du\right\}dW(s).  \label{eq:3.6}
\end{equation}
\end{prop}

\begin{proof}
For $-\infty<a<b<\infty$ with $(a,b]\subset I$, we have 
\begin{equation*}
X(b)-X(a)=\int_{-\infty}^{\infty}
\left\{\int_{I}I_{(a,b]}(u)c(u-s)du\right\}dW(s),
\end{equation*}
which implies (\ref{eq:3.6}) for $f=I_{(a,b]}$. The general case follows
easily from this.
\end{proof}

\begin{prop}
\label{prop:3.3} Let $f\in \mathcal{H}_I(X)$ such that $f\geq 0$, and let $%
f_{n} \ (n=1,2,\dots)$ be a sequence of simple functions on $I$ such that $%
0\leq f_{n}\uparrow f$ a.e. Then, in $M(X)$, 
\begin{equation*}
\lim_{n\to\infty}\int_{-\infty}^{\infty}f_{n}(s) dX(s)
=\int_{-\infty}^{\infty}\left\{ \int_{I}f(u)c(u-s)du\right\}dW(s).
\end{equation*}
\end{prop}

\begin{proof}
By Proposition \ref{prop:3.2} and the monotone convergence theorem, we have 
\begin{equation*}
\begin{aligned} &\left\Vert \int_{I}f_{n}(s) dX(s) -\int_{-\infty}^{\infty}
\left\{\int_{I}f(u)c(u-s) du\right\} dW(s) \right\Vert^2 \\ &\le
\int_{-\infty}^{\infty} \left\{\int_{I}(f(u) - f_{n}(u))c(u-s)du\right\}^2ds
\ \downarrow\ 0, \qquad n\to\infty. \end{aligned}
\end{equation*}
Thus the proposition follows.
\end{proof}

For a real-valued function $f$ on $I$, we write $f(x)=f^+(x)-f^-(x)$, where 
\begin{equation*}
f^+(x):=\max(f(x),0),\quad f^-(x):=\max(-f(x),0), \qquad x\in I.
\end{equation*}

\begin{defn}
\label{def:3.4} For $f\in H_{I}(X)$, we define 
\begin{equation*}
\int_{I}f(s)dX(s):=\lim_{n\rightarrow \infty
}\int_{I}f_{n}^{+}(s)dX(s)-\lim_{n\rightarrow \infty
}\int_{I}f_{n}^{-}(s)dX(s)\quad \mbox{in $M(X)$},
\end{equation*}%
where $\{f_{n}^{+}\}$ and $\{f_{n}^{-}\}$ are arbitrary sequences of
non-negative simple functions on $I$ such that $f_{n}^{+}\uparrow f^{+}$, $%
f_{n}^{-}\uparrow f^{-}$, as $n\rightarrow \infty $, a.e. 
\end{defn}

From the definition above, we see that $\int_If(s)dX(s)\in M_I(X)$ for $f\in%
\mathcal{H}_I(X)$. The next proposition follows immediately from Proposition %
\ref{prop:3.3}.

\begin{prop}
\label{prop:3.5} The equality\/ $(\ref{eq:3.6})$ also holds for $f\in%
\mathcal{H}_I(X)$.
\end{prop}

\subsection{Infinite past prediction formulas}

We denote by $\mathcal{D}(\mathbf{R})$ the space of all $\phi \in C^{\infty
}(\mathbf{R})$ with compact support, endowed with the usual topology. For a
random distribution $Y$ (cf.\ \cite[\S 2]{I2} and \cite[\S 2]{AIK}%
), we write $DY$ for its derivative. For $t\in\mathbf{R}$, we write $%
M_{(-\infty,t]}(Y)$ for the closed linear hull of $\{Y(\phi ):\phi \in%
\mathcal{D}(\mathbf{R}),\ \mathrm{supp}\ \phi \subset (-\infty,t]\}$ in $%
L^2(\Omega,\mathcal{F},P)$. Notice that $M_I(X)$ here coincides with that
defined above.

As in \cite[Proposition 2.4]{IA}, we have the next proposition.

\begin{prop}
\label{prop:3.6} The derivative\/ $DX$ of\/ $(X(t))$ is a purely
nondeterministic stationary random distribution, and $(W(t):t\in\mathbf{R})$
is a canonical Brownian motion of $DX$ in the sense that $%
M_{(-\infty,t]}(DX)=M_{(-\infty,t]}(DW)$ for every $t\in\mathbf{R}$.
\end{prop}

See \S 5 for the proof.

Here is the infinite past prediction formula for $\int_{t}^{\infty}f(s)dX(s)$%
.

\begin{thm}
\label{thm:3.7} For $t\in [0,\infty)$ and $f\in\mathcal{H}_{[t,\infty)}(X)$,
the following assertions hold:

\begin{itemize}
\item[\textrm{(a)}] $\int_{0}^{\infty}b(t-\cdot,\tau)f(t+\tau)d\tau\in 
\mathcal{H}_{(-\infty,t]}(X)$.

\item[\textrm{(b)}] $P_{(-\infty,t]}\int_{t}^{\infty}f(s)dX(s)
=\int_{-\infty}^{t}\left\{\int_{0}^{\infty}b(t-s,\tau)f(t+\tau)d\tau
\right\}dX(s)$.
\end{itemize}
\end{thm}

\begin{proof}
Since $f\in\mathcal{H}_{[t,\infty)}(X)$ iff $\vert f\vert\in\mathcal{H}%
_{[t,\infty)}(X)$, we may assume $f\ge 0$. Since 
\begin{equation}
c(u)=0,\qquad t\le 0,  \label{eq:3.7}
\end{equation}
it follows from (\ref{eq:3.5}) and the Fubini--Tonelli theorem that, for $%
s<t $, 
\begin{equation}
\begin{aligned} \int_{t}^{\infty}f(u)c(u-s)du &=\int_{0}^{\infty}d\tau
f(t+\tau)\int_{0}^{t-s}c(t-s-u)b(u,\tau)du\\
&=\int_{-\infty}^{t}duc(u-s)\int_{0}^{\infty}b(t-u,\tau)f(t+\tau)d\tau.
\end{aligned}  \label{eq:3.8}
\end{equation}
Thus we obtain (a). By Proposition \ref{prop:3.6} and \cite[Proposition 2.3\
(2)]{AIK}, we have 
\begin{equation}
M_{(-\infty,t]}(X)=M_{(-\infty,t]}(DW).  \label{eq:3.9}
\end{equation}
This and Proposition \ref{prop:3.5} yield 
\begin{equation*}
P_{(-\infty,t]}\int_{t}^{\infty}f(s)dX(s)
=\int_{-\infty}^{t}\left\{\int_{t}^{\infty} f(u)c(u-s)du\right\}dW(s).
\end{equation*}
By (\ref{eq:3.7}), (\ref{eq:3.8}) and Proposition \ref{prop:3.5}, the
integral on the right-hand side is 
\begin{equation*}
\begin{aligned} &\int_{-\infty}^{t}\left\{\int_{-\infty}^{t}duc(u-s)
\int_{0}^{\infty}b(t-u,\tau)f(t+\tau)d\tau\right\}dW(s)\\
&=\int_{-\infty}^{t}\left\{\int_{0}^{\infty}b(t-s,\tau)f(t+\tau)d\tau\right%
\} dX(s). \end{aligned}
\end{equation*}
Thus (b) follows.
\end{proof}

By putting $f(s)=I_{(t_1,T]}(s)$ in Theorem \ref{thm:3.7} (b), we
immediately obtain the next infinite past prediction formula for $(X(t))$.

\begin{thm}
\label{thm:3.8} Let $0\le t_1< T<\infty$. Then $\int_{0}^{T-t_1}b(t_1-\cdot,%
\tau)d\tau\in \mathcal{H}_{(-\infty,t_1]}(X)$ and the infinite past
prediction formula $(\ref{eq:1.11})$ holds.
\end{thm}

Using the Hilbert space isomorphism $\theta :M(X)\rightarrow M(X)$
characterized by $\theta (X(t))=X(-t)$ for $t\in \mathbf{R}$, we obtain the
next theorem from Theorem \ref{thm:3.7} (see the proof of \cite[Theorem 3.6]%
{AIK}).

\begin{thm}
\label{thm:3.9} For $t\in [0,\infty)$ and $f\in\mathcal{H}_{[t,\infty)}(X)$,
the following assertions hold:

\begin{itemize}
\item[\textrm{(a)}] $\int_{0}^{\infty}b(t+\cdot,\tau)f(t+\tau)d\tau\in 
\mathcal{H}_{[-t,\infty)}(X)$.

\item[\textrm{(b)}] $P_{[-t,\infty)}\int_{-\infty}^{-t}f(-s)dX(s)
=\int_{-t}^{\infty}\left\{\int_{0}^{\infty}b(t+s,\tau)f(t+\tau)d\tau
\right\}dX(s)$.
\end{itemize}
\end{thm}

As in \cite[Definition 2.2]{AIK}, we define another Brownian motion $%
(W^*(t):t\in\mathbf{R})$ by 
\begin{equation}
W^*(t):=\theta(W(-t)), \qquad t\in\mathbf{R}.  \label{eq:3.10}
\end{equation}

\begin{prop}
\label{prop:3.10} Let $I$ be a closed interval of\/ $\mathbf{R}$ and let $%
f\in\mathcal{H}_{I}(X)$. Then 
\begin{equation*}
\int_{I}f(s)dX(s)
=\int_{-\infty}^{\infty}\left\{\int_{I}f(u)c(s-u)du\right\}dW^*(s).
\end{equation*}
\end{prop}

The proof of Proposition \ref{prop:3.10} is the same as that of \cite[%
Proposition 3.5]{AIK}, whence we omit it. We need Theorem \ref{thm:3.9} and
Proposition \ref{prop:3.10} in the next section.

\begin{exmp}
\label{exmp:3.11} As in Example \ref{exmp:2.3}, we consider fBm $(B_{H}(t))$
with $1/2<H<1$. Then the MA$(\infty )$ coefficient $c(t)$ is given by 
\begin{equation}
c(t)=t^{H-(3/2)}\frac{1}{\Gamma (H-\frac{1}{2})},\qquad t>0,  \label{eq:3.11}
\end{equation}%
so that $\int_{0}^{\infty }e^{izt}c(t)dt=(-iz)^{(1/2)-H}$ for $\Im z>0$. From (%
\ref{eq:3.2}), we have 
\begin{equation*}
\int_{0}^{\infty }e^{izt}\alpha (t)dt=(-iz)^{H-3/2}.
\end{equation*}%
Hence, $\alpha (t)=t^{(1/2)-H}/\Gamma (\frac{3}{2}-H)$, so that the AR$%
(\infty )$ coefficient $a(t)$ is given by 
\begin{equation}
a(t)=t^{-(H+\frac{1}{2})}\frac{H-\frac{1}{2}}{\Gamma (\frac{3}{2}-H)},\qquad
t>0.  \label{eq:3.12}
\end{equation}%
By the change of variable $u=sv$, $%
\int_{0}^{s}(s-u)^{H-(3/2)}(t+u)^{-H-(1/2)}du$ becomes 
\begin{equation*}
s^{H-\frac{1}{2}}t^{-H-\frac{1}{2}}\int_{0}^{1}(1-v)^{H-\frac{3}{2}%
}\{1+(s/t)v\}^{-H-\frac{1}{2}}dv=\frac{1}{(H-\frac{1}{2})}\left( \frac{s}{t}%
\right) ^{H-\frac{1}{2}}\frac{1}{t+s},
\end{equation*}%
where we have used the equality 
\begin{equation*}
\int_{0}^{1}(1-v)^{p-1}(1+xv)^{-p-1}dv=\frac{1}{p(x+1)},\qquad p>0,\ x>-1.
\end{equation*}%
Thus 
\begin{equation}
b(t,s)=\frac{\sin \{\pi (H-\frac{1}{2})\}}{\pi }\left( \frac{s}{t}\right)
^{H-\frac{1}{2}}\frac{1}{t+s},\qquad t>0,\ s>0;  \label{eq:3.13}
\end{equation}%
and so, from Theorem \ref{thm:3.8}, we see that, for $0\leq t<T$, 
\begin{equation*}
\begin{aligned} &E\left[B_H(T) \vert\ \sigma(B_H(s): -\infty< s\le
t)\right]\\ &=B_H(t)+\frac{\sin\{\pi (H-\frac{1}{2})\}}{\pi}
\int_{-\infty}^{t}\left\{\int_{0}^{T-t}
\left(\frac{\tau}{t-s}\right)^{H-\frac{1}{2}}\frac{1}{t-s+\tau}d\tau\right%
\}dB_H(s). \end{aligned}
\end{equation*}%
This prediction formula was obtained in \cite[Theorem 3.1]{GN} by a
different method. 
\end{exmp}


\section{Finite past prediction problems}
\label{sec:4}

In this section, we assume (\ref{eq:1.1})--(\ref{eq:1.7}) and (\ref{eq:1.10}%
). Notice that (\ref{eq:1.6}) with (\ref{eq:1.7}) implies (\ref{eq:3.1}) as
well as (\ref{eq:2.3}), whence (\ref{eq:2.1}). For $t_0$, $t_1$, and $T$ in (%
\ref{eq:1.10}), we put 
\begin{equation*}
t_2:=t_0+t_1,\quad t_3:=T-t_1.
\end{equation*}

\subsection{Alternating projections to the past and future}

For $n\in \mathbf{N}$, we define the orthogonal projection operator $P_{n}$
by 
\begin{equation*}
P_{n}:=%
\begin{cases}
P_{(-\infty ,t_{1}]}, & n=1,3,5,\dots , \\ 
P_{[-t_{0},\infty )}, & n=2,4,6,\dots .%
\end{cases}%
\end{equation*}%
It should be noted that $\{P_{n}\}_{n=1}^{\infty }$ is merely an alternating
sequence of projection operators, first to $M_{(-\infty ,t_{1}]}(X)$, then
to $M_{[-t_{0},\infty )}(X)$, and so on. This sequence plays a key role in
the proof of the finite past prediction formula for $(X(t))$.

For $t,s\in (0,\infty)$ and $n\in\mathbf{N}$, we define $%
b_n(t,s)=b_n(t,s;t_2)$ iteratively by 
\begin{equation}
\begin{cases}
b_1(t,s):=b(t,s), &  \\ 
b_n(t,s):=\int_{0}^{\infty}b(t,u)b_{n-1}(t_2+u,s)du, & n=2,3,\dots.%
\end{cases}
\label{eq:4.1}
\end{equation}

\begin{prop}
\label{prop:4.1} For $f\in\mathcal{H}_{[t_1,\infty)}(X)$, the following
assertions hold:

\begin{itemize}
\item[\textrm{(a)}] $\int_{0}^{\infty}b_n(t_1-\cdot,\tau)f(t_1+\tau)d\tau \in%
\mathcal{H}_{(-\infty,t_1]}(X)$ for $n=1,3,5,\dots$.

\item[\textrm{(b)}] $\int_{0}^{\infty}b_n(t_0+\cdot,\tau)f(t_1+\tau)d\tau \in%
\mathcal{H}_{[-t_0,\infty)}(X)$ for $n=2,4,6,\dots$.
\end{itemize}
\end{prop}

\begin{proof}
We may assume that $f\ge 0$. By Theorem \ref{thm:3.7}, (a) holds for $n=1$.
By the Fubini--Tonelli theorem, we have, for $s>-t_0$, 
\begin{equation*}
\int_{0}^{\infty}dub(t_0+s,u)\int_{0}^{\infty}b_1(t_2+u,\tau)
f(t_1+\tau)d\tau =\int_{0}^{\infty}b_2(t_0+s,\tau)f(t_1+\tau)d\tau.
\end{equation*}
Hence, by Theorem \ref{thm:3.9}, we have (b) for $n=2$. Repeating this
procedure, we obtain the proposition.
\end{proof}

Let $f\in \mathcal{H}_{[t_1,\infty)}(X)$. By Proposition \ref{prop:4.1}, we
may define the random variables $G_n(f)$ by 
\begin{equation*}
G_n(f):= 
\begin{cases}
\int_{-t_0}^{t_1}\left\{
\int_{0}^{\infty}b_n(t_1-s,\tau)f(t_1+\tau)d\tau\right\}dX(s), & n=1,3,\dots,
\\ 
\int_{-t_0}^{t_1}\left\{
\int_{0}^{\infty}b_n(t_0+s,\tau)f(t_1+\tau)d\tau\right\}dX(s), & n=2,4,\dots.%
\end{cases}%
\end{equation*}
We may also define the random variables $\epsilon_n(f)$ by $%
\epsilon_0\left(f\right):=\int_{t_1}^{\infty}f(s)dX(s)$ and 
\begin{equation*}
\epsilon_n(f):= 
\begin{cases}
\int_{-\infty}^{-t_0}\left\{
\int_{0}^{\infty}b_n(t_1-s,\tau)f(t_1+\tau)d\tau\right\}dX(s), & n=1,3,\dots,
\\ 
\int_{t_1}^{\infty}\left\{
\int_{0}^{\infty}b_n(t_0+s,\tau)f(t_1+\tau)d\tau\right\}dX(s), & n=2,4,\dots.%
\end{cases}%
\end{equation*}

\begin{prop}
\label{prop:4.2} Let $f\in\mathcal{H}_{[t_1,\infty)}(X)$ and $n\in\mathbf{N}$%
. Then 
\begin{equation}
P_nP_{n-1}\cdots P_1\int_{t_1}^{\infty}f(s)dX(s)=\epsilon_n\left(f\right)
+\sum_{k=1}^{n}G_k(f).  \label{eq:4.2}
\end{equation}
\end{prop}

We can prove (\ref{eq:4.2}) using Proposition \ref{prop:4.1} and the facts 
\begin{align}
&M_{[-t_0,t_1]}(X)\subset M_{(-\infty,t_1]}(X) \cap M_{[-t_0,\infty)}(X),
\label{eq:4.3} \\
&G_k\in M_{[-t_0,t_1]}(X),\qquad k=1,2,\dots.  \label{eq:4.4}
\end{align}
Since the proof is similar to that of \cite[Proposition 4.4]{AIK}, we omit
the details.

We are about to investigate the limit of (\ref{eq:4.2}) as $n\to\infty$ (see
Lemma \ref{lem:4.9} below).

For $f\in\mathcal{H}_{[t_1,\infty)}(X)$ and $s>0$, we define $%
D_n(s,f)=D_n(s,f;t_1,t_2)$ by 
\begin{equation*}
D_n(s,f):= 
\begin{cases}
\int_{0}^{\infty}c(u)f(t_1+s+u)du, & n=0, \\ 
\int_{0}^{\infty}duc(u) \int_{0}^{\infty}b_n(t_2+u+s,\tau)f(t_1+\tau) d\tau,
& n=1,2,\dots.%
\end{cases}%
\end{equation*}
From the proof of the next proposition, we see that these integrals converge
absolutely. Recall $(W^*(t))$ from (\ref{eq:3.10}).

\begin{prop}
\label{prop:4.3} Let $f\in\mathcal{H}_{[t_1,\infty)}(X)$. Then 
\begin{equation*}
P_{n+1}^{\bot}\epsilon_n\left(f\right) = 
\begin{cases}
\int_{t_1}^{\infty}D_n(s-t_1,f)dW(s), & n=0,2,4,\dots, \\ 
\int_{-\infty}^{-t_0}D_n(-t_0-s,f)dW^*(s), & n=1,3,5,\dots.%
\end{cases}%
\end{equation*}
\end{prop}

\begin{proof}
By (\ref{eq:3.9}) and Proposition \ref{prop:3.5}, 
\begin{equation*}
P_1^{\bot}\epsilon_0\left(f\right)
=\int_{t_1}^{\infty}\left\{\int_{s}^{\infty}f(u)c(u-s)du\right\}dW(s)
=\int_{t_1}^{\infty}D_0(s-t_1,f)dW(s).
\end{equation*}
Thus the assertion holds for $n=0$. Let $n=1,3,\dots$. Then, by Proposition %
\ref{prop:3.10}, 
\begin{equation*}
\epsilon_n\left(f\right)
=\int_{-\infty}^{\infty}\left\{\int_{-\infty}^{-t_0}duc(s-u)
\int_{0}^{\infty}b_n(t_1-u,\tau)f(t_1+\tau)d\tau\right\}dW^*(s).
\end{equation*}
Hence, using \cite[Proposition 2.3\ (7)]{AIK} and (\ref{eq:3.7}), 
\begin{equation*}
\begin{aligned} P_{n+1}^{\bot}\epsilon_n\left(f\right)
&=\int_{-\infty}^{-t_0}\left\{\int_{-\infty}^{s}duc(s-u)
\int_{0}^{\infty}b_n(t_1-u,\tau)f(t_1+\tau)d\tau\right\}dW^*(s)\\
&=\int_{-\infty}^{-t_0}\left\{\int_{0}^{\infty}duc(u)
\int_{0}^{\infty}b_n(t_2+u-t_0-s,\tau)f(t_1+\tau)d\tau\right\}dW^*(s)\\
&=\int_{-\infty}^{-t_0}D_n(-t_0-s,f)dW^*(s). \end{aligned}
\end{equation*}
Thus we obtain the assertion for $n=1,3,\dots$. The proof for $n=2,4,\dots $
is similar; and so we omit it.
\end{proof}

From Propositions \ref{prop:4.2} and \ref{prop:4.3}, we immediately obtain
the next proposition (cf.~the proof of \cite[Proposition 4.9]{AIK}).

\begin{prop}
\label{prop:4.4} Let $f\in\mathcal{H}_{[t_1,\infty)}(X)$. Then the following
assertions hold:

\begin{itemize}
\item[\textrm{(a)}] $\Vert P_{1}^{\bot}\int_{t_1}^{\infty}f(s)dX(s)\Vert^2
=\int_{0}^{\infty}D_0(s,f)^2ds$.

\item[\textrm{(b)}] $\Vert P^{\bot}_{n+1}P_nP_{n-1}\cdots P_1
\int_{t_1}^{\infty}f(s)dY(s)\Vert^2 =\int_{0}^{\infty}D_n(s,f)^2ds$ for $%
n=1,2,\dots$.
\end{itemize}
\end{prop}

We write $Q$ for the orthogonal projection operator from $M(X)$ onto the
intersection $M_{(-\infty,t_1]}(X)\cap M_{[-t_0,\infty)}(X)$. Then, by von
Neumann's alternating projection theorem (see, e.g., \cite[Theorem 9.20]{P}), 
we have $Q=\slim_{n\to\infty}P_nP_{n-1}\cdots P_1$%
. Using this, (\ref{eq:4.3}) and Proposition \ref{prop:4.4}, we immediately
obtain the next proposition (cf.\ the proof of \cite[Proposition 4.9\ (3)]%
{AIK}).

\begin{prop}
\label{prop:4.5} Let $f\in\mathcal{H}_{[t_1,\infty)}(X)$. Then $%
\lim_{n\to\infty}\int_{0}^{\infty}D_n(s,f)^2ds=0$.
\end{prop}

We need the next proposition.

\begin{prop}
\label{prop:4.6} Let $f\in\mathcal{H}_{[t_1,\infty)}(X)$. Then, for $t>0$
and $n=0,1,\dots$, we have 
\begin{equation*}
\int_{0}^{\infty}b_{n+1}(t,\tau)f(t_1+\tau)d\tau
=\int_{0}^{\infty}a(t+u)D_n(u,f)du.
\end{equation*}
\end{prop}

\begin{proof}
We may assume $f\ge 0$. By the Fubini--Tonelli theorem, we have, for $t>0$, 
\begin{equation*}
\begin{aligned} &\int_{0}^{\infty}b_1(t,\tau)f(t_1+\tau)d\tau
=\int_{0}^{\infty}\left\{
\int_{0}^{\tau}c(\tau-u)a(t+u)du\right\}f(t_1+\tau)d\tau\\
&\quad=\int_{0}^{\infty} a(t+u)
\left\{\int_{0}^{\infty}c(\tau)f(t_1+u+\tau)d\tau\right\}du
=\int_{0}^{\infty}a(t+u)D_0(u,f)du. \end{aligned}
\end{equation*}
Thus the assertion holds for $n=0$. Now we assume that $n\ge 1$. Since we
have 
\begin{equation*}
b_{n+1}(t,\tau)=\int_{0}^{\infty}a(t+v)\left\{\int_{0}^{%
\infty}c(u)b_n(t_2+u+v,\tau)du\right\}dv, \quad t, \tau>0,
\end{equation*}
we obtain the assertion, again using the Fubini--Tonelli theorem.
\end{proof}

For $t,s>0$, we define $k(t,s)=k(t,s;t_2)$ by 
\begin{equation*}
k(t,s):=\int_{0}^{\infty}c(t+u)a(t_2+u+s)du.
\end{equation*}
Notice that $k(t,s)<\infty$ for $t,s>0$ since $k(t,s)\le
c(t)\int_{t_2+s}^{\infty}a(u)du$.

\begin{prop}
\label{prop:4.7} Let $f\in\mathcal{H}_{[t_1,\infty)}(X)$. Then 
\begin{equation*}
P_{n+1}\epsilon_n\left(f\right) = 
\begin{cases}
\int_{-\infty}^{t_1}\left\{\int_{0}^{\infty}k(t_1-s,u)
D_{n-1}(u,f)du\right\}dW(s), & n=2,4,\dots, \\ 
\int_{-t_0}^{\infty}\left\{\int_{0}^{\infty}k(t_0+s,u)D_{n-1}(u,f)du\right\}
dW^*(s), & n=1,3,\dots.%
\end{cases}%
\end{equation*}
\end{prop}

\begin{proof}
We assume $n=2,4,\dots$. Then, by Propositions \ref{prop:3.5} and \ref%
{prop:4.6}, we have 
\begin{equation*}
\begin{aligned} P_{n+1}\epsilon_n\left(f\right)
&=\int_{-\infty}^{t_1}\left\{\int_{t_1}^{\infty}duc(u-s)
\int_{0}^{\infty}b_n(t_0+u,\tau)f(t_1+\tau)d\tau\right\}dW(s)\\
&=\int_{-\infty}^{t_1}\left\{\int_{0}^{\infty}dv c(t_1-s+v)
\int_{0}^{\infty}a(t_2+v+u)D_{n-1}(u,f)du\right\}dW(s)\\
&=\int_{-\infty}^{t_1}\left\{\int_{0}^{\infty}k(t_1-s,u)
D_{n-1}(u,f)du\right\}dW(s). \end{aligned}
\end{equation*}
The proof of the case $n=1,3,\dots$ is similar.
\end{proof}

We need the next $L^2$-boundedness theorem.

\begin{thm}
\label{thm:4.8} Let $p\in (0,1/2)$ and let $\ell(\cdot)$ be a slowly varying
function at infinity. Let $C(\cdot)$ and $A(\cdot)$ be nonnegative and
decreasing functions on $(0,\infty)$. We assume $C(\cdot)\in L^1_{\mathrm{loc%
}}[0,\infty)$ and $A(0+)<\infty$. We also assume 
\begin{align*}
&A(t)\sim t^{-(1+p)}\ell(t)p, \qquad t\to\infty, \\
&C(t)\sim \frac{t^{-(1-p)}}{\ell(t)}\cdot \frac{\sin(p\pi)}{\pi}, \qquad
t\to\infty,
\end{align*}
and put $K(x,y):=\int_{0}^{\infty}C(x+u)A(u+y)du$ for $x, y>0$. 
Then
\[
\sup_{x>0}\int_{0}^{\infty}K(x,y)\left(x/y\right)^{1/2}dy <\infty, \qquad 
\sup_{y>0}\int_{0}^{\infty}K(x,y)\left(y/x\right)^{1/2}dx <\infty.
\]
In particular, 
the integral operator $K$ defined by $(Kf)(x):=\int_{0}^{\infty}K(x,y)f(y)dy$
for $x>0$ is a bounded operator on $L^2((0,\infty),dy)$.
\end{thm}

We omit the proof of Theorem \ref{thm:4.8} which is similar to that of 
\cite[Theorem 5.1]{IA}.

By putting $z=iy$ in (\ref{eq:3.2}), we get 
\begin{equation*}
y\left\{\int_{0}^{\infty}e^{-yt}c(t)dt\right\}
\left\{\int_{0}^{\infty}e^{-yt}\alpha(t)dt\right\}=1, \qquad y>0.
\end{equation*}
By Karamata's Tauberian theorem (cf.\ \cite[Theorem 1.7.6]{BGT}) applied to
this, (\ref{eq:2.3}) implies 
$\alpha(t)\sim t^{-(H-\frac{1}{2})}/\{\ell(t)\Gamma((3/2)-H)\}$ as 
$t\to\infty$. This and the monotone density theorem give 
\begin{equation}
a(t)\sim \frac{t^{-(H+\frac{1}{2})}}{\ell(t)} \cdot\frac{(H-\frac{1}{2})}{%
\Gamma(\frac{3}{2}-H)}, \qquad t\to\infty.  \label{eq:4.5}
\end{equation}

The next lemma is a key to our arguments.

\begin{lem}
\label{lem:4.9} Let $f\in\mathcal{H}_{[t_1,\infty)}(X)$. Then $\Vert
\epsilon_n(f)\Vert\to 0$ as $n\to\infty$.
\end{lem}

\begin{proof}
It follows from (\ref{eq:2.3}), (\ref{eq:4.5}) and Theorem \ref{thm:4.8}
below that the integral operator $K$ defined by $Kf(t):=\int_{0}^{%
\infty}k(t,s)f(s)ds$ is a bounded operator on $L^2((0,\infty),ds)$. Hence,
by Propositions \ref{prop:4.3}, \ref{prop:4.5} and \ref{prop:4.7}, we have 
\begin{equation*}
\begin{aligned} \left\Vert \epsilon_n\left(f\right)\right\Vert^2
&=\int_{0}^{\infty}D_n(s,f)^2ds+\int_{0}^{\infty}
\left\{\int_{0}^{\infty}k(s,u)D_{n-1}(u,f)du\right\}^2ds\\ &\le
\int_{0}^{\infty}D_{n}(s,f)^2ds +\Vert
K\Vert^2\int_{0}^{\infty}D_{n-1}(s,f)^2ds \to 0, \qquad n\to\infty.
\end{aligned}
\end{equation*}
Thus the lemma follows.
\end{proof}

We can now state the conclusions of the arguments above.

\begin{thm}
\label{thm:4.10} The following assertions hold:

\begin{itemize}
\item[\textrm{(a)}] $M_{[-t_0,t_1]}(X)= M_{(-\infty,t_1]}(X)\cap
M_{[-t_0,\infty)}(X)$.

\item[\textrm{(b)}] $P_{[-t_0,t_1]} =\slim_{n\to\infty}P_nP_{n-1}\cdots P_1$.

\item[\textrm{(c)}] $\Vert P_{[-t_0,t_1]}^{\bot}Z\Vert^2 =\left\Vert
P_1^{\bot}Z\right\Vert^2 +\sum_{n=1}^{\infty}\left\Vert
(P_{n+1})^{\bot}P_n\cdots P_1Z\right\Vert^2$ for $Z\in M(X)$.
\end{itemize}
\end{thm}

We can prove Theorem \ref{thm:4.10} using Proposition \ref{prop:4.2} and
Lemma \ref{lem:4.9}. Since the proof is similar to that of \cite[Theorem 4.6]%
{AIK}, we omit the details.

\subsection{Finite past prediction formulas}

We define $h(s,u)=h(s,u;t_2)$ by 
\begin{equation}
h(s,u):=\sum_{k=1}^{\infty}\left\{b_{2k-1}(t_2-s,u) +b_{2k}(s,u)\right\},
\qquad 0<s<t_2,\ u>0.  \label{eq:4.6}
\end{equation}

Here is the finite past prediction formula for $\int_{t_1}^{\infty}f(s)dX(s)$%
.

\begin{thm}
\label{thm:4.11} Let $f\in\mathcal{H}_{[t_1,\infty)}(X)$. Then the following
assertions hold:

\begin{itemize}
\item[\textrm{(a)}] $\int_{0}^{\infty}h(t_0+\cdot,u)f(t_1+u)du \in\mathcal{H}%
_{[-t_0,t_1]}(X)$.

\item[\textrm{(b)}] $P_{[-t_0,t_1]}\int_{t_1}^{\infty}f(s)dX(s)
=\int_{-t_0}^{t_1}\left\{\int_{0}^{\infty}h(t_0+s,u)f(t_1+u)du\right\}dX(s)$.

\item[\textrm{(c)}] $\Vert P_{[-t_0,t_1]}^{\bot}
\int_{t_1}^{\infty}f(s)dX(s)\Vert^2
=\sum_{n=0}^{\infty}\int_{0}^{\infty}D_n(s,f)^2ds$.
\end{itemize}
\end{thm}

\begin{proof}
We may assume that $f\ge 0$. By Theorem \ref{thm:4.10} (b), Proposition \ref%
{prop:4.2} and Lemma \ref{lem:4.9}, we have, in $M(X)$, 
\begin{equation*}
\begin{aligned} &P_{[-t_0,t_1]}\int_{t_1}^{\infty}f(s)dX(s)
=\lim_{n\to\infty}P_nP_{n-1}\cdots P_1\int_{t_1}^{\infty}f(s)dX(s)\\
&=\lim_{n\to\infty}
\int_{-t_0}^{t_1}\left\{\int_{0}^{\infty}h_n(t_0+u,v)f(t_1+v)dv\right%
\}dX(s), \end{aligned}
\end{equation*}
where, for $0<s<t_2$ and $u>0$, we define $h_n(s,u)=h_n(s,u;t_2)$ by 
\begin{equation*}
h_n(s,u)= 
\begin{cases}
b_1(t_2-s,u)+b_2(s,u)+\cdots +b_n(t_2-s,u), & n=1,3,5,\dots, \\ 
b_1(t_2-s,u)+b_2(s,u)+\cdots +b_n(s,u), & n=2,4,6,\dots.%
\end{cases}%
\end{equation*}
Since $h_n(s,u)\uparrow h(s,u)$, we obtain (a) and (b) using the monotone 
convergence theorem. Finally, (c) follows from
Theorem \ref{thm:4.10} (c) and Proposition \ref{prop:4.4}.
\end{proof}

For $s,u>0$, we define $D_n(s)=D_n(s;t_2,t_3)$ by 
\begin{equation*}
D_n(s) :=\int_{0}^{\infty}du c(u)\int_{0}^{t_3}b_n(t_2+u+s,\tau)d\tau,
\qquad n=1,2,\dots.
\end{equation*}

Here are the solutions to the finite past prediction problems for $(X(t))$.

\begin{thm}
\label{thm:4.12} The finite past prediction formula $(\ref{eq:1.12})$ and
the following equality for the mean-square prediction error hold: 
\begin{equation*}
\left\Vert P_{[-t_0,t_1]}^{\bot}X(T)\right\Vert^2
=\int_{0}^{T-t_1}g(s)^2ds+\sum_{n=1}^{\infty}\int_{0}^{\infty}D_n(s)^2ds.
\end{equation*}
\end{thm}

\begin{proof}
We put $f(s)=I_{(t_{1},T]}(s)$. Then $\int_{t_{1}}^{\infty
}f(s)dX(s)=X(T)-X(t_{1})$ and 
\begin{equation*}
\int_{0}^{\infty
}h(t_{0}+s,u)f(t_{1}+u)du=\int_{0}^{t_{3}}h(t_{0}+s,u)du,\qquad
-t_{0}<s<t_{1}.
\end{equation*}%
We also have $D_{n}(s,f)=D_{n}(s)$ for $n=1,2,\dots $ and $%
D_{0}(s,f)=g(t_{3}-s)$. Thus the theorem follows from Theorem \ref{thm:4.11}.
\end{proof}


\section{AR$(\infty)$-type equations}
\label{sec:5}

In this section, we consider the AR$(\infty)$-type equations for 
$(X(t))$ in (\ref{eq:1.1}) and $(\tilde{X}(t))$ in (\ref{eq:1.13}). 
For a Borel measure $\tau$ on $(0,\infty)$ satisfying 
$\int_0^{\infty}(1+s)^{-1}\tau(ds)<\infty$, we write
\[
F_{\tau}(z):=\int_0^{\infty}\frac{1}{\lambda-iz}\tau(d\lambda),\qquad \Im z\ge 0.
\]

First, we consider the process $X=(X(t))$ in (\ref{eq:1.1}) with (\ref{eq:1.2})--(\ref{eq:1.5}), 
$(\ref{eq:2.1})$ and $(\ref{eq:3.1})$. 
Let $f_t(s):=g(t-s)-g(-s)=\int_{-s}^{t-s}c(u)du$ for $t,s\in\mathbf{R}$.

\begin{lem}\label{lem:5.1}
Let $t\in\mathbf{R}$. 
Then the Fourier transform of $f_{t}(\cdot)$ in the $L^{2}$-sense 
is equal to $(i\xi)^{-1}(1-e^{-it\xi})F_{\nu}(\xi)$:
\begin{equation}
\frac{(1-e^{-it\xi})}{i\xi}F_{\nu}(\xi) 
=\underset{M\to\infty}{\mathrm{l.i.m.}}\int_{-M}^{M}
e^{-is\xi}f_{t}(s) ds.
\label{eq:5.1}
\end{equation}
\end{lem}

\begin{proof}
Since $\int_{-\infty}^{\infty}\vert f_t(s)\vert^2ds<\infty$, 
the limit on the right-hand side of (\ref{eq:5.1}) exists. Therefore, it is enough to justify 
the following point-wise convergence:
\begin{equation}
\frac{\left( 1-e^{-it\xi }\right) }{i\xi }F_{{}\nu }\left( \xi \right) =%
\lim_{M\to \infty}\int_{-M}^{M}e^{-is\xi}f_{t}(s)ds,\qquad \xi \neq 0. 
\label{eq:5.2}
\end{equation}
Now, if $-M\le t\le M$, then 
\begin{align*}
&\int_{-M}^{M}e^{-is\xi }f_t(s)ds 
=\int_{-M}^{M}dse^{-is\xi}\int_0^tc(u-s)du
=\int_0^tdu\int_{-M}^{M}e^{-is\xi}c(u-s)ds\\
&\qquad\qquad\qquad=\int_{0}^{t}due^{-iu\xi}\int_{u-M}^{u+M}e^{iv\xi }c(v) dv
=\int_{0}^{t}due^{-iu\xi}\int_{0}^{u+M}e^{iv\xi }c(v) dv
\end{align*}
because $u-M\leq 0\le u+M$ for $u$ between $0$ and $t$, and $c(s)=0$ for $s\leq 0$. 
However, 
\begin{align*}
&\int_{0}^{t}due^{-iu\xi}\int_{0}^{u+M}e^{is\xi}c(s)ds
=\int_{0}^{t}due^{-iu\xi}
\int_{0}^{\infty}\frac{1-e^{(i\xi - \lambda)(u + M)}}{\lambda -i\xi}\nu(d\lambda)  \\
&\qquad\qquad\qquad\qquad\qquad 
=\frac{(1-e^{-it\xi})}{i\xi}F_{\nu}(\xi) 
-e^{i\xi M}\int_0^t du\int_{0}^{\infty}\frac{e^{-\lambda(u+M)}}{\lambda -i\xi}\nu(d\lambda),
\end{align*}
so that, for $\xi \neq 0$, 
\[
\left| \frac{(1-e^{-it\xi})}{i\xi}F_{\nu}(\xi) -\int_{-M}^{M}e^{-is\xi }f_{t}(s) ds\right| 
\le t\int_{0}^{\infty }\frac{e^{-\lambda M}}{|\lambda -i\xi| 
}\nu(d\lambda)\ \downarrow\ 0,\quad M\rightarrow \infty .
\]
Thus, (\ref{eq:5.2}) holds.
\end{proof}

For the Brownian motion $W=(W(t))$ in (\ref{eq:1.1}), 
let $DW(\phi)=\int_{-\infty }^{\infty }\hat{\phi}(\xi
)Z_{DW}(d\xi )$ with $\phi \in \mathcal{D}(\mathbf{R})$ be the spectral
decomposition of $DW$ as a stationary random distribution, where 
$\hat{\phi}(\xi):=\int_{-\infty }^{\infty }e^{-it\xi }\phi (\xi)d\xi$ 
and $Z_{DW}$ 
is the associated complex-valued random measure such that 
$E[Z_{DW}(A)\overline{Z_{DW}(B)}]=(2\pi)^{-1}\int_{A\cap B}d\xi$ 
(see It\^o \cite{It}). By Lemma \ref{lem:5.1} and the Parseval-type 
formula for the homogeneous random measure $Z_{DW}$, we obtain 
$X(t)=\int_{-\infty}^{\infty}[(1-e^{-it\xi})/(i\xi)]F_{\nu}(\xi)Z_{DW}(d\xi)$, 
whence
\begin{equation}
DX(\phi)=\int_{-\infty}^{\infty}\hat{\phi}(\xi)F_{\nu}(\xi)Z_{DW}(d\xi),\qquad 
\phi \in \mathcal{D}(\mathbf{R}).
\label{eq:5.3}
\end{equation}
Let $\rho_{DX}$ be the spectral measure of $DX$: 
$E[X(\phi)\overline{X(\psi)}]
=\int_{-\infty}^{\infty}\hat{\phi}(\xi)\overline{\hat{\psi}}(\xi)\rho_{DX}(d\xi)$. 
Then, from (\ref{eq:5.3}), we see that 
$\rho_{DX}(d\xi)=(2\pi)^{-1}\vert F_{\nu}(\xi)\vert^2d\xi$. 
Thus, $DX$ has the spectral density 
$\Delta_{DX}(\xi):=(2\pi)^{-1}\vert F_{\nu}(\xi)\vert^2$. 
Since, for $z=x+iy$ with $y>0$, we have
\[
\Re\left\{ F_{\nu}(z)\right\} =\int_{0}^{\infty }\frac{s+y}{(s+y)^{2}+x^{2}}\nu(ds)>0,
\]
the function $F_{\nu}(z)$ is an outer function on the upper half plane $\Im z>0$:
\begin{equation}
F_{\nu}(z)=\exp\left\{\frac{1}{\pi i}\int_{-\infty}^{\infty}
\frac{1+\xi z}{\xi-z}\cdot\frac{\log \vert F_{\nu}(\xi)\vert}{1+\xi^2}d\xi\right\},
\qquad \Im z>0.
\label{eq:5.4}
\end{equation}
In particular, Proposition \ref{prop:3.6} follows from this and (\ref{eq:5.3}).

We also have the next lemma.

\begin{lem}\label{lem:5.2}
The following equality holds:
\[
\int_{0}^{\infty}e^{izt}c(t)dt
=\sqrt{2\pi}\exp\left\{\frac{1}{2\pi i}\int_{-\infty}^{\infty}
\frac{1+\xi z}{\xi-z}\cdot\frac{\log \vert \Delta_{DX}(\xi)\vert}{1+\xi^2}d\xi\right\},
\qquad \Im z>0.
\]
\end{lem}

\begin{proof}
Since $F_{\nu}(z)=\int_{0}^{\infty}e^{izt}c(t)dt$ and 
$\vert F_{\nu}(\xi)\vert=\{2\pi\Delta_{DX}(\xi)\}^{1/2}$, 
the lemma follows from (\ref{eq:5.4}).
\end{proof}

From Lemma \ref{lem:5.2}, we see that the kernel $c(\cdot)$ is uniquely determined 
by $DX$, whence $(X(t))$, as claimed in \S \ref{subsec:3.2}

Let $D^2X:=D(DX)$. For the AR$(\infty)$ kernel $\alpha(\cdot)$ in \S \ref{subsec:3.1}, 
we define the convolution $\alpha\ast D^2X$, which is also
a stationary random distribution, by 
\begin{equation}
(\alpha\ast D^2X)(\phi )
:=\underset{M\to\infty}{\mathrm{l.i.m.}}
\int_{0}^{M}\alpha(u)D^2X(\tau_{u}\phi )du,\qquad
\phi \in \mathcal{D}(\mathbf{R}),
\label{eq:5.5}
\end{equation}
where $\tau_{u}\phi(t):=\phi (t+u)$ and the integral on the right-hand
side is an $M(X)$-valued Bochner integral. Then, by \cite[Proposition 2.3]{I2} 
and (\ref{eq:5.3}), we have
\[
(\alpha\ast D^2X)(\phi)
=-\int_{-\infty}^{\infty}i\xi F_{\mu}(\xi)F_{\nu}(\xi)\hat{\phi}(\xi)Z_{DW}(d\xi).
\]
However, since (\ref{eq:3.2}) implies $-i\xi F_{\mu}(\xi)F_{\nu}(\xi)=1$ for 
$\xi\ne 0$, we see that $X$ satisfies
\begin{equation}
\alpha\ast D^2X=DW.
\label{eq:5.6}
\end{equation}
More precisely, we have the next theorem.

\begin{thm}\label{thm:5.3}
The process $(X(t))$ is the only stationary-increment process with $X(0)=0$ satisfying the following two conditions:
\begin{enumerate}
\item the stationary random distribution $DX$ is purely nondeterministic;
\item $(X(t))$ satisfies $(\ref{eq:5.6})$.
\end{enumerate}
\end{thm}

The proof of Theorem \ref{thm:5.3} is similar to that of \cite[Theorem 2.6]{AI2}, whence we omit it. 
Notice that (\ref{eq:5.6}) can be written formally as the following 
AR$(\infty)$-type equation:
\begin{equation}
\int_{-\infty}^t\alpha(t-s)\frac{d^2X}{ds^2}(s)ds=\frac{dW}{dt}(t).
\label{eq:5.7}
\end{equation}

\begin{exmp}
\label{exmp:5.4} 
Let $(B_H(t))$ be the fBm in (\ref{eq:2.5}) with $1/2<H<1$. Then, by Example \ref{exmp:3.11}, 
we have $\alpha(t)=t^{(1/2)-H}/\Gamma(\frac{3}{2}-H)$ for $t>0$, whence (\ref{eq:5.7}) becomes
\[
\frac{1}{\Gamma(\frac{3}{2}-H)}\int_{-\infty}^{t}\frac{1}{(t-s)^{H-\frac{1}{2}}}
\cdot\frac{d^2B_H}{ds^2}(s)ds=\frac{dW}{dt}(t).
\]
\end{exmp}

Next, we turn to $\tilde{X}=(\tilde{X}(t))$ in (\ref{eq:1.13}) with (\ref{eq:1.14}). We assume 
that $\tilde{\nu}$ is a Borel measure on $(0,\infty)$ satisfying the following conditions:
\[
\int_0^{\infty}\frac{1}{1+s}\tilde{\nu}(ds)<\infty,\quad 
\tilde{\nu}((0,\infty))=\int_{0}^{\infty}\frac{1}{s}\tilde{\nu}(ds)=\infty,\quad 
\int_0^1\tilde{c}(t)^2dt<\infty.
\]
By \cite[Theorem 3.2]{IA}, there exists a unique Borel measure $\tilde{\mu}$ on $(0,\infty)$ 
satisfying
\[
\int_0^{\infty}\frac{1}{1+s}\tilde{\mu}(ds)<\infty,\quad 
\tilde{\mu}((0,\infty))=\int_{0}^{\infty}\frac{1}{s}\tilde{\mu}(ds)=\infty,
\]
and $-izF_{\tilde{\nu}}(z)F_{\tilde{\mu}}(z)=1$ for $\Im z>0$. 
If we define
\[
\tilde{\alpha}(t):=\int_0^{\infty}e^{-st}\tilde{\mu}(ds),\qquad t>0,
\]
then the last equality becomes
\begin{equation}
-iz\left\{ \int_{0}^{\infty }e^{izt}\tilde{c}(t)dt\right\} \left\{ \int_{0}^{\infty
}e^{izt}\tilde{\alpha} (t)dt\right\} =1,\qquad \Im z>0.
\label{eq:5.8}
\end{equation}

By \cite[(2.3)]{IA}, we have
\begin{equation}
D\tilde{X}(\phi)=\int_{-\infty}^{\infty}\hat{\phi}(\xi)(-i\xi)F_{\tilde{\nu}}(\xi)Z_{DW}(d\xi),\qquad 
\phi \in \mathcal{D}(\mathbf{R}),
\label{eq:5.9}
\end{equation}
whence, in the same way as the proof of \cite[Proposition 5.1]{I1}, we get
\[
(\tilde{\alpha}\ast D\tilde{X})(\phi)
=-\int_{-\infty}^{\infty}i\xi F_{\tilde{\mu}}(\xi)F_{\tilde{\nu}}(\xi)\hat{\phi}(\xi)Z_{DW}(d\xi),\qquad 
\phi \in \mathcal{D}(\mathbf{R}),
\]
where the convolution $\tilde{\alpha}\ast D\tilde{X}$ is defined in the same way as (\ref{eq:5.5}). 
However, since $-i\xi F_{\tilde{\mu}}(\xi)F_{\tilde{\nu}}(\xi)=1$ for 
$\xi\ne 0$, we see that $(\tilde{X}(t))$ satisfies
\begin{equation}
\tilde{\alpha}\ast D\tilde{X}=DW.
\label{eq:5.10}
\end{equation}
Notice that the equation (\ref{eq:5.10}) can be written formally as the 
following 
AR$(\infty)$-type equation:
\begin{equation}
\int_{-\infty}^t\tilde{\alpha}(t-s)\frac{d\tilde{X}}{ds}(s)ds=\frac{dW}{dt}(t).
\label{eq:5.11}
\end{equation}
We can also prove an analogue of Theorem \ref{thm:5.3} for $(\tilde{X}(t))$, 
which we omit in this paper.

\begin{exmp}
\label{exmp:5.5} 
Let $(B_H(t))$ be the fBm in (\ref{eq:2.5}) with $0<H<1/2$. Then, by 
\cite[Example 3.9]{IA}, 
we have $\tilde{\alpha}(t)=t^{-(1/2)-H}/\Gamma(\frac{1}{2}-H)$ for $t>0$, whence 
(\ref{eq:5.11}) becomes
\[
\frac{1}{\Gamma(\frac{1}{2}-H)}\int_{-\infty}^{t}\frac{1}{(t-s)^{H+\frac{1}{2}}}
\cdot\frac{dB_H}{ds}(s)ds=\frac{dW}{dt}(t).
\]
\end{exmp}

\ 

{\bf Acknowledgements.}\quad 
We would like to express our gratitude to an anonymous referee for useful suggestions.

\


\end{document}